\documentclass{article}
\usepackage{amsmath, amssymb}
\newtheorem{theorem}{Theorem}

\newtheorem{definition}[theorem]{Definition}

\newtheorem{lemma}[theorem]{Lemma}
\newtheorem{proposition}[theorem]{Proposition}
\newtheorem{remark}[theorem]{Remark}
\newenvironment{proof}[1][Proof]{\noindent\textbf{#1.} }{\ \rule{0.5em}{0.5em}}

\begin{document}
\title{Two-Scale limit of the solution to a Robin Problem in Perforated Media
 \footnote{\hspace*{-7mm} \textit{2000 Mathematics Subject Classification:}
35B27.
\newline
 \textit{Keywords:} Homogenization, Two scale
convergence, Robin boundary condition. }}
\author{Abdelhamid Ainouz\\
 Department of Mathematics, University of Sciences \\ and Technology
Houari Boumediene, \\ Po Box 32 El Alia Bab Ezzouar 16111 Algiers,
Algeria. \\aainouz@usthb.dz}

\maketitle

\begin{abstract}
The two scale convergence of the solution to a Robin's type-like
problem of a stationary diffusion problem in a periodically
perforated domain is investigated. It is shown that the Robin's
problem converges to a problem associated to a new operator which is
the sum of a standard homogenized operator plus an extra first order
"strange" term; its appearance is due to the non-symmetry of the
diffusion matrix and to the non rescaled resistivity.
\end{abstract}

\section{Introduction}
Periodic homogenization in perforated media with Robin boundary
conditions prescribed on the boundary of the holes has been
extensively studied by many
authors and we refer for instance to\ \cite{bpc}, \cite{ciodon1}, \cite%
{ciodon2}, \cite{pas1},\ .... In this paper we study the
stationary diffusion equation in a periodic perforated body where
the heat flow is proportional to the temperature field on the
boundary of the holes with a resistivity having zero average value
on the boundary of the reference hole. In \cite{bpc}, the authors
studied a model problem for a second-order symmetric elliptic
operator in a periodically perforated domain with the Robin
boundary condition prescribed on the boundary of the holes. They
use the asymptotic expansion technique \cite{blp}, \cite{san},
\cite{tar} to obtain the homogenized  problem and  they construct
correctors to justify the expansion and then estimate the error.
In this paper, we consider a problem but with another
configuration, namely the holes may or may not be connected and
the boundary of the holes may intersect the exterior boundary of
the body. Moreover we assume that the matrix diffusion of the
second-order operator may or may not be symmetric. We use the
two-scale convergence technique \cite{all}, \cite{lnw},
\cite{ngue}, \cite{ngue1} to obtain the two-scale limit system.
After the decoupling technique, we show that the homogenized
problem contains a convective term. Its appearance is due
essentially to the general character of the matrix diffusion and
on the fact that the resistivity function is not rescaled as
usually assumed when dealing with two-scale convergence on
periodic surfaces, see for instance \cite{ain}, \cite{ain1},
\cite{adh}, \cite{nr},...

The paper is organized as follows: in section \ref{1}, we define
the geometry of the perforated body and we give the Robin
boundary-value problem setting. Section \ref{2} is aimed at
showing the existence and unicity of the solution of the Robin
problem and obtaining a priori estimates of the solution. The
asymptotic limit via two-scale convergence procedure is analyzed
in section \ref{3}. We obtain the homogenized boundary-value
problem which is a second order elliptic operator containing first
and zero order terms. The latter term is classical see e.g.
\cite{adh}. The first order term is a convection one and it is
null when the diffusion matrix is symmetric or with constant
coefficients.

\section{Setting of the Problem\label{1}}

Let $\Omega $ be a bounded domain in $\mathbb{R}^{n}$ of variable
$x=\left( x_{i}\right) _{1\leq i\leq n}\ $($n\geq 2$) with a
smooth boundary $\Gamma $ and $\varepsilon $ a real parameter
taking values in a sequence of positive numbers tending to zero.
As usual in periodic homogenization let $Y=\left[ 0,1\right] ^{n}$
be the generic unit cell of periodicity in the auxiliary space
$\mathbb{R}^{n}$ of variable $y=\left( y_{i}\right) _{1\leq i\leq
n}$.
The cell $Y$ is identified to the unit torus $\mathbb{R}^{n}/\mathbb{Z}^{n}$%
. A function defined on $\mathbb{R}^{n}$ is said to be
$Y$-periodic if it is periodic of period $1$ in each $y_{i}$
variable with $1\leq i\leq n$. In the sequel we will suppose that
any function defined on $Y$ is extended periodically to the whole
space $\mathbb{R}^{n}$. If $E\left( Z\right) $ is a function space
(where $Z$ is a subset of $Y$) we denote $E_{\#}\left( Z\right)
:=\left\{ w\in E\left( Z\right) \text{; }w\text{ is extended
periodically to }\mathbb{R}^{n}\right\} $.

Let $H$, the reference hole, be an open subset of $Y$ with a
smooth boundary $\Sigma $ and set $Y_{s}=Y\backslash cl(H)$ where
$cl\left( \cdot \right) $ denotes the closure. Thus $Y$ is
partitioned as $Y=Y_{s}\cup \Sigma \cup H$. Note that we do not
require that $H$ is strictly included in $Y$. As a consequence the
periodic extension of $H$ may or may not be connected. Let
us denote $\chi \left( y\right) $ the characteristic function of $Y_{s}$ in $%
Y$. We define the perforated material
\begin{equation*}
\Omega _{\varepsilon }=\left\{ x\in \Omega ;\chi \left(
\frac{x}{\varepsilon }\right) =1\right\}
\end{equation*}%
and the one codimensional periodic surface
\begin{equation*}
\Sigma _{\varepsilon }=\left\{ x\in \Omega ;\frac{x}{\varepsilon
}\in \Sigma \right\} .
\end{equation*}

Here $\Omega _{\varepsilon }$ represents the matrix or the solid part of $%
\Omega ,$ by opposition to the holes or the void part that is
represented by the open subset $H_{\varepsilon }:=\Omega
\backslash cl\left( \Omega _{\varepsilon }\right) $. By
construction, all of these holes are identical and they are
periodically distributed in $\Omega $ with period $\varepsilon $
in each $x_{i}$-direction. Since we use the two-scale convergence
method we do not require that the boundary $\Sigma _{\varepsilon
}=\partial H_{\varepsilon }$ does not intersect $\Gamma $. As in
\cite{all}, We shall use the natural extension by zero of any
function defined on $\Omega _{\varepsilon }$.

Let $f_{\varepsilon }$ be a given function in $L^{2}\left( \Omega
_{\varepsilon }\right) $, $g_{\varepsilon }$ be given in
$L^{2}\left( \Sigma _{\varepsilon }\right) $ such that
\begin{equation}
\Vert f_{\varepsilon }\Vert _{L^{2}\left( \Omega _{\varepsilon }\right) }+%
\sqrt{\varepsilon }\Vert g_{\varepsilon }\Vert _{L^{2}\left(
\Sigma _{\varepsilon }\right) }\leq C\text{.}  \label{dc}
\end{equation}%
Here and throughout this paper $C$ denotes a positive constant
independent of $\varepsilon $.

Let $A\left( x,y\right) =\left( a_{ij}\right) _{1\leq i,j\leq n}$
be a real-valued matrix function defined on $\Omega \times Y$,
$Y$-periodic in the second variable $y$ such that there exists two
positive constants $m$, $M
$ independent of $\varepsilon $ satisfying the following inequality:%
\begin{equation}
m\mid \zeta \mid ^{2}\leq (A\zeta ,\zeta )\leq M\mid \zeta \mid
^{2} \label{ac}
\end{equation}%
for all $\zeta \in \mathbb{R}^{n}$. We suppose that the matrix $A$ lies in $%
C\left( \Omega ;L_{\#}^{\infty }\left( Y\right) \right) ^{n^{2}}$.
We note that no symmetry condition on $A$ is assumed. Let $\mu
\left( y\right) \in L_{\#}^{\infty }\left( Y_{s}\right) $ such
that $\int_{Y_{s}}\mu \left( y\right) \geq \mu _{0}>0$ where $\mu
_{0}$ is independent of $\varepsilon $.
Let $\alpha $ be a $Y$-periodic measurable bounded function defined on $%
\Sigma $ such that
\begin{equation}
\int_{\Sigma }\alpha \left( y\right) d\sigma \left( y\right)
=0\text{.} \label{as1}
\end{equation}%
Let us decompose the function $\alpha $ into its positive and
negative parts
as follows:%
\begin{equation*}
\alpha =\alpha ^{+}-\alpha ^{-},\ \alpha ^{+}=\max \left( \alpha
,0\right) \text{, }\alpha ^{-}=\max \left( -\alpha ,0\right)
\text{.}
\end{equation*}%
Assume that the positive part of $\alpha $ satisfies the
condition:
\begin{equation}
\alpha ^{+}\left( y\right) \geq \alpha _{0}>0\text{ a.e. in
}\Sigma \text{.} \label{ca}
\end{equation}%
Let us consider the following Robin boundary value problem:

\begin{eqnarray}
-div\left( A_{\varepsilon }\nabla u_{\varepsilon }\right) +\mu
_{\varepsilon }u_{\varepsilon } &=&f_{\varepsilon }\text{ in
}\Omega _{\varepsilon }, \label{eq1} \\ \left( A_{\varepsilon
}\nabla u_{\varepsilon }\right) \cdot \nu _{\varepsilon }+\alpha
_{\varepsilon }u_{\varepsilon } &=&\varepsilon g_{\varepsilon
}\text{ on }\Sigma _{\varepsilon },  \label{eq2} \\ u_{\varepsilon
} &=&0\text{ on }\Gamma   \label{eq3}
\end{eqnarray}%
where
\begin{equation*}
A_{\varepsilon }\left( x\right) =A\left( x,\frac{x}{\varepsilon
}\right) ,\ \mu _{\varepsilon }\left( x\right) =\mu \left(
\frac{x}{\varepsilon }\right) ,\ \alpha _{\varepsilon }\left(
x\right) =\alpha \left( \frac{x}{\varepsilon }\right)
\end{equation*}%
and $\nu _{\varepsilon }$ is the unit outward normal to $\Omega
_{\varepsilon }$.

This problem can be regarded as a simplified model of the
condensation of stream in a periodic cooling structure (see
\cite{adh}). We can also consider as a model for treatment
planning hyperthemia in microvascular tissue- see, e.g. \cite{dh}.
The boundary condition (\ref{eq2}) means that the heat flow
$\left( A_{\varepsilon }\nabla u_{\varepsilon }\right) \cdot \nu
_{\varepsilon }$ is proportional to the temperature
$u_{\varepsilon }$ with a periodic resistivity given by the
function $\alpha _{\varepsilon }$. In many situations, the
resistivity function is taken to be $\varepsilon ^{m}\alpha
_{\varepsilon }$. Since the operator is of order $2$ the
interesting cases are then $m=-2,-1,0,1$ and $2$. The case $m=2$
is trivial since we obtain the classical homogenized equation.This
can be seen easily by using the asymptotic expansion method. The
case $m=1$  with $\alpha _{\varepsilon }\geq 0$ has been studied
by \cite{adh}, \cite{nr} using the two-scale convergence
technique. In this situation $\alpha _{\varepsilon }$ is rescaled
since the surface $\Sigma _{\varepsilon }$ is of codimension $1$ \
Here we study the case $m=0$, i.e. a non rescaled resistivity. We
use the same technique but with $\alpha _{\varepsilon }$ changing
sign. We show that the assumptions (\ref{as1}), (\ref{ca}) and a
non-symmetric $A_{\varepsilon }\left( x\right) $ contribute to the
description of the effective thermal conductivity with convection.
The case $m=-1$ will be studied in a forthcoming paper. \ Note
that the case $m=-2$ is also trivial since it yields that the
effective thermal conductivity is $0$.

\section{Study of the Problem and A priori Estimates\label{2}}

Let
\begin{equation*}
V_{\varepsilon }=\left\{ v\in H^{1}\left( \Omega _{\varepsilon
}\right) \text{; }v=0\text{ on }\Gamma \right\}
\end{equation*}%
equipped with the scalar product
\begin{equation*}
(u,v)_{V_{\varepsilon }}=\int_{\Omega _{\varepsilon }}\nabla
u\left( x\right) \nabla v\left( x\right) dx
\end{equation*}%
and the associated norm $\Vert u\Vert _{V_{\varepsilon
}}=(u,u)_{V_{\varepsilon }}^{1/2}$ which is equivalent to the
$H^{1}$-norm thanks to the Poincar\'{e} inequality$.$ The
variational formulation of the
boundary-value problem (\ref{eq1})-(\ref{eq3}) reads as follows:%
\begin{equation}
\left\{
\begin{array}{c}
\text{For each }\varepsilon >0\text{, find }u_{\varepsilon }\in
V_{\varepsilon }\text{ such that } \\ a_{\varepsilon }\left(
u_{\varepsilon },v\right) =L_{\varepsilon }\left(
v\right) \text{ for any }v\in V_{\varepsilon },%
\end{array}%
\right.  \label{wf}
\end{equation}%
where $a_{\varepsilon }\left( \cdot ,\cdot \right) $ is the
bilinear form defined on $V_{\varepsilon }\times V_{\varepsilon }$
by:
\begin{eqnarray*}
a_{\varepsilon }\left( u,v\right) &=&\int_{\Omega _{\varepsilon
}}A_{\varepsilon }(x)\nabla u\left( x\right) \nabla v\left(
x\right) dx+\int_{\Omega _{\varepsilon }}\mu _{\varepsilon }\left(
x\right) u\left( x\right) v\left( x\right) dx \\ &&+\int_{\Sigma
_{\varepsilon }}\alpha _{\varepsilon }\left( x\right) u\left(
x\right) v\left( x\right) d\sigma _{\varepsilon }\left( x\right)
\end{eqnarray*}%
and $L_{\varepsilon }\left( \cdot \right) $ is the linear form defined on $%
V_{\varepsilon }$ by:
\begin{equation*}
L_{\varepsilon }\left( v\right) =\int_{\Omega _{\varepsilon
}}f_{\varepsilon }(x)v\left( x\right) dx+\varepsilon \int_{\Sigma
_{\varepsilon }}g_{\varepsilon }\left( x\right) v\left( x\right)
d\sigma _{\varepsilon }\left( x\right) .
\end{equation*}

\begin{lemma}
\label{lcs}There exists a positive constant $C_{s}$\ independent of $%
\varepsilon $ such that for every $v\in V_{\varepsilon }$ and for every $%
\delta >0$ we have
\begin{equation}
\Vert v\Vert _{L^{2}\left( \Sigma _{\varepsilon }\right) }^{2}\leq
C_{s} \left[ \left( \delta \varepsilon \right) ^{-1}\Vert v\Vert
_{L^{2}\left( \Omega _{\varepsilon }\right) }^{2}+\left( \delta
\varepsilon \right) \Vert \nabla v\Vert _{\left( L^{2}\left(
\Omega _{\varepsilon }\right) \right) ^{n}}^{2}\right] .
\label{cs}
\end{equation}
\end{lemma}

\begin{proof}
Let us introduce the notation
\begin{equation*}
v_{\varepsilon }^{k}\left( x\right) =v\left( \varepsilon \left(
k+y\right) \right)
\end{equation*}%
where $k\in K_{\varepsilon }=\left\{ k\in
\mathbb{Z}^{n};\varepsilon \left(
Y+k\right) \cap \Omega \neq \emptyset \right\} $. By the change of variable $%
x=\varepsilon \left( k+y\right) $ we have
\begin{equation*}
\int_{\Sigma _{\varepsilon }}v^{2}\left( x\right) d\sigma
_{\varepsilon }\left( x\right) =\underset{k\in K_{\varepsilon
}}{\sum }\int_{\varepsilon \left( \Sigma +k\right) }v^{2}\left(
x\right) d\sigma _{\varepsilon }\left(
x\right) =\varepsilon ^{n-1}\underset{k\in K_{\varepsilon }}{\sum }%
\int_{\Sigma }\left[ v_{\varepsilon }^{k}\left( y\right) \right]
^{2}d\sigma \left( y\right) .
\end{equation*}%
From the trace theorem we see that for every $\delta >0$
\begin{eqnarray*}
\int_{\Sigma }\left[ v_{\varepsilon }^{k}\left( y\right) \right]
^{2}d\sigma \left( y\right) &\leq &C_{s}\left[ \delta
^{-1}\int_{Y_{s}+k}\left[ v_{\varepsilon }^{k}\left( y\right)
\right] ^{2}dy+\delta \int_{Y_{s}+k}|\nabla _{y}v_{\varepsilon
}^{k}\left( y\right) |^{2}dy\right]
\\
&\leq &\frac{C_{s}}{\varepsilon ^{n}}\left[ \delta
^{-1}\int_{\varepsilon \left( Y_{s}+k\right) }v\left( x\right)
^{2}dx+\delta \varepsilon ^{2}\int_{\varepsilon \left(
Y_{s}+k\right) }|\nabla _{x}v\left( x\right) |^{2}dx\right] .
\end{eqnarray*}%
Hence

\begin{eqnarray*}
\int_{\Sigma _{\varepsilon }}v^{2}\left( x\right) d\sigma
_{\varepsilon
}\left( x\right)  &\leq &\varepsilon ^{n-1}\frac{C_{s}}{\varepsilon ^{n}}%
[\delta ^{-1}\underset{k\in K_{\varepsilon }}{\sum
}\int_{\varepsilon \left( Y_{s}+k\right) }v\left( x\right) ^{2}dx+
\\
&&\delta \varepsilon ^{2}\underset{k\in K_{\varepsilon }}{\sum }%
\int_{\varepsilon \left( Y_{s}+k\right) }|\nabla _{x}v\left(
x\right) |^{2}dx] \\ &\leq &C_{s}\left[ \left( \delta \varepsilon
\right) ^{-1}\int_{\Omega _{\varepsilon }}v^{2}\left( x\right)
dx+\left( \delta \varepsilon \right) \int_{\Omega _{\varepsilon
}}|\nabla v\left( x\right) |^{2}dx\right] \text{.}
\end{eqnarray*}%
The Lemma is proved.
\end{proof}

\begin{lemma}
Let $\sqrt{\mu _{0}m}>C_{s}\Vert \alpha \Vert _{L^{\infty }\left(
\Sigma
\right) }$ where $C_{s}$ is the constant given in Lemma \ref{lcs}. Then $%
a_{\varepsilon }\left( \cdot ,\cdot \right) $ is coercive on
$V_{\varepsilon }$.
\end{lemma}

\begin{proof}
Let $v\in V_{\varepsilon }$. Then using (\ref{ac}), we have
\begin{equation*}
a_{\varepsilon }\left( v,v\right) \geq m\int_{\Omega _{\varepsilon
}}|\nabla v\left( x\right) |^{2}dx+\mu _{0}\int_{\Omega
_{\varepsilon }}v\left( x\right) ^{2}dx-\Vert \alpha \Vert
_{L^{\infty }\left( \Sigma \right) }\int_{\Sigma _{\varepsilon
}}v\left( x\right) ^{2}d\sigma _{\varepsilon }\left( x\right) .
\end{equation*}%
By (\ref{cs}) we see that for every $\delta >0$
\begin{eqnarray*}
a_{\varepsilon }\left( v,v\right)  &\geq &\left( m-\delta
\varepsilon C_{s}\Vert \alpha \Vert _{L^{\infty }\left( \Sigma
\right) }\right) \int_{\Omega _{\varepsilon }}|\nabla v\left(
x\right) |^{2}dx+ \\ &&\left( \mu _{0}-\left( \delta \varepsilon
\right) ^{-1}C_{s}\Vert \alpha \Vert _{L^{\infty }\left( \Sigma
\right) }\right) \int_{\Omega _{\varepsilon }}\left[ v\left(
x\right) \right] ^{2}dx.
\end{eqnarray*}

Choosing $\delta =\dfrac{1}{\varepsilon }\sqrt{\dfrac{m}{\mu
_{0}}}$. Then
\begin{equation*}
a_{\varepsilon }\left( v,v\right) \geq c_{0}\left[ \Vert \nabla
v\Vert _{\left( L^{2}\left( \Omega _{\varepsilon }\right) \right)
^{n}}^{2}+\Vert v\Vert _{L^{2}\left( \Omega _{\varepsilon }\right)
}^{2}\right]
\end{equation*}%
where $c_{0}$ is the positive constant given by
\begin{equation*}
c_{0}=\left( 1-\frac{C_{s}\Vert \alpha \Vert _{L^{\infty }\left(
\Sigma \right) }}{\sqrt{m\mu _{0}}}\right) \min \left( m,\mu
_{0}\right) >0.
\end{equation*}%
This completes the proof.
\end{proof}

In the sequel we shall assume that the condition $\sqrt{\mu _{0}m}%
>C_{s}\Vert \alpha \Vert _{L^{\infty }\left( \Sigma \right) }$ is fulfilled.

\begin{proposition}
The variational formulation (\ref{wf}) admits a unique solution $%
u_{\varepsilon }\in V_{\varepsilon }$. Moreover we have the a
priori
estimates.%
\begin{equation}
\Vert \nabla u_{\varepsilon }\Vert _{\left( L^{2}\left( \Omega
_{\varepsilon }\right) \right) ^{n}}+\Vert u_{\varepsilon }\Vert
_{L^{2}\left( \Omega _{\varepsilon }\right) }+\Vert u_{\varepsilon
}\Vert _{L^{2}\left( \Sigma _{\varepsilon }\right) }\leq C.
\label{ape}
\end{equation}
\end{proposition}

\begin{proof}
The existence and uniqueness is a straightforward application of Lemma \ref%
{lcs} and the Lax-Milgram Lemma. It remains to prove the a priori estimates (%
\ref{ape}). Take $v=$\bigskip $u_{\varepsilon }$ in (\ref{wf}). We have%
\begin{eqnarray*}
&&\int_{\Omega _{\varepsilon }}\left( A_{\varepsilon }\nabla
u_{\varepsilon }\nabla u_{\varepsilon }+\mu _{\varepsilon
}u_{\varepsilon }^{2}\right) dx+\int_{\Sigma _{\varepsilon
}}\alpha _{\varepsilon }^{+}u_{\varepsilon }^{2}d\sigma
_{\varepsilon }\left( x\right)  \\ &=&\int_{\Omega _{\varepsilon
}}f_{\varepsilon }u_{\varepsilon }dx+\int_{\Sigma _{\varepsilon
}}\left( \varepsilon g_{\varepsilon }+\alpha _{\varepsilon
}^{-}u_{\varepsilon }\right) u_{\varepsilon }d\sigma _{\varepsilon
}\left( x\right) .
\end{eqnarray*}%
Let us denote
\begin{equation*}
A_{\varepsilon }\left( u_{\varepsilon }\right) :=\Vert \nabla
u_{\varepsilon }\Vert _{\left( L^{2}\left( \Omega _{\varepsilon
}\right) \right) ^{n}}^{2}+\Vert u_{\varepsilon }\Vert
_{L^{2}\left( \Omega _{\varepsilon }\right) }^{2}+\Vert
u_{\varepsilon }\Vert _{L^{2}\left( \Sigma _{\varepsilon }\right)
}^{2}.
\end{equation*}%
Then using (\ref{ac}) and (\ref{ca}) we obtain
\begin{equation}
A_{\varepsilon }\left( u_{\varepsilon }\right) \leq \frac{1}{c_{1}}%
(\int_{\Omega _{\varepsilon }}f_{\varepsilon }u_{\varepsilon
}dx+\int_{\Sigma _{\varepsilon }}\left( \varepsilon g_{\varepsilon
}+\alpha _{\varepsilon }^{-}u_{\varepsilon }\right) u_{\varepsilon
}d\sigma _{\varepsilon }\left( x\right) )  \label{ione}
\end{equation}%
where $c_{1}=\min \left( m,\mu _{0},\alpha _{0}\right) >0$.
Applying Young's inequality on the right hand side of
(\ref{ione}), we get
\begin{eqnarray}
A_{\varepsilon }\left( u_{\varepsilon }\right)  &\leq &\frac{1}{c_{1}}[\frac{%
\beta ^{2}}{2}\Vert f_{\varepsilon }\Vert _{L^{2}\left( \Omega
_{\varepsilon }\right) }^{2}+\frac{1}{2\beta ^{2}}\Vert
u_{\varepsilon }\Vert _{L^{2}\left( \Omega _{\varepsilon }\right)
}^{2}  \notag \\ &&+\frac{\gamma ^{2}\varepsilon ^{2}}{2}\Vert
g_{\varepsilon }\Vert _{L^{2}\left( \Omega _{\varepsilon }\right)
}^{2}+\left( \frac{1}{2\gamma ^{2}}+\Vert \alpha \Vert _{L^{\infty
}\left( \Sigma \right) }\right) \Vert u_{\varepsilon }\Vert
_{L^{2}\left( \Sigma _{\varepsilon }\right) }^{2}]. \label{itwo}
\end{eqnarray}%
But in view of (\ref{cs}) inequality (\ref{itwo}) becomes%
\begin{eqnarray*}
A_{\varepsilon }\left( u_{\varepsilon }\right)  &\leq &\frac{1}{c_{1}}%
[\left( \frac{1}{2\beta ^{2}}+\left( \frac{\varepsilon }{2\gamma
^{2}}+\Vert
\alpha \Vert _{L^{\infty }\left( \Sigma \right) }\right) \frac{C_{s}}{%
\varepsilon \delta }\right) \Vert u_{\varepsilon }\Vert
_{L^{2}\left( \Omega _{\varepsilon }\right) }^{2} \\ &&+\left(
\frac{\varepsilon }{2\gamma ^{2}}+\Vert \alpha \Vert _{L^{\infty
}\left( \Sigma \right) }\right) C_{s}\varepsilon \delta \Vert
\nabla u_{\varepsilon }\Vert _{\left( L^{2}\left( \Omega
_{\varepsilon }\right) \right) ^{n}}^{2}]+C.
\end{eqnarray*}%
Now, appropriate choice of $\beta ,\gamma ,\delta $ yields
\begin{equation*}
A_{\varepsilon }\left( u_{\varepsilon }\right) =\Vert \nabla
u_{\varepsilon }\Vert _{\left( L^{2}\left( \Omega _{\varepsilon
}\right) \right) ^{n}}^{2}+\Vert u_{\varepsilon }\Vert
_{L^{2}\left( \Omega _{\varepsilon }\right) }^{2}+\Vert
u_{\varepsilon }\Vert _{L^{2}\left( \Sigma _{\varepsilon }\right)
}^{2}\leq C.
\end{equation*}%
The Proposition is now proved.
\end{proof}

One is led to determine the homogenized problem of
(\ref{eq1})-(\ref{eq3}). Namely we study the limiting behavior of
the solutions $u_{\varepsilon }$ as $\varepsilon $ tends to zero.
This the subject of the next section.

\section{ Homogenization Procedure\label{3}}

We shall use the well-known two-scale convergence method that we
briefly recall here the definition and the main results.

\subsection{Two-scale Convergence}

\begin{definition}
\begin{enumerate}
\item A sequence $v_{\varepsilon }$ in $L^{2}(\Omega )$ \emph{two-scale}
converges to $v_{0}(x,y)\in L^{2}((\Omega \times Y)$ and we denote this $%
v_{\varepsilon }\rightrightarrows v_{0}$ if for any $\varphi
(x,y)\in L^{2}\left( \Omega ;C_{\#}\left( Y\right) \right) $,
\begin{equation*}
\lim_{\varepsilon \rightarrow 0}\int_{\Omega }v_{\varepsilon
}(x)\varphi \left( x,\frac{x}{\varepsilon }\right) dx=\int_{\Omega
}\int_{Y}v_{0}(x,y)\varphi (x,y)dydx.
\end{equation*}

\item A sequence $v_{\varepsilon }$ in $L^{2}(\Sigma _{\varepsilon })$ \emph{%
two-scale} converges to $v_{0}(x,y)\in L^{2}((\Omega \times \Sigma
)$ and we denote this $v_{\varepsilon
}\overset{S}{\rightrightarrows }v_{0}$ if for any $\varphi
(x,y)\in C\left( \overline{\Omega };C_{\#}\left( Y\right) \right)
$,
\begin{equation*}
\lim_{\varepsilon \rightarrow 0}\int_{\Sigma _{\varepsilon
}}\varepsilon v_{\varepsilon }(x)\varphi \left(
x,\frac{x}{\varepsilon }\right) d\sigma _{\varepsilon }\left(
x\right) =\int_{\Omega }\int_{\Sigma }v_{0}(x,y)\varphi
(x,y)d\sigma _{\varepsilon }\left( y\right) dx.
\end{equation*}
\end{enumerate}
\end{definition}

\begin{proposition}
\label{p1}

\begin{enumerate}
\item For any uniformly bounded sequence $v_{\varepsilon }$ in $L^{2}\left(
\Omega \right) $ one can extract a subsequent still denoted by
$\varepsilon $
and a two-scale limit $v_{0}\in L^{2}((\Omega \times Y)$ such that $%
v_{\varepsilon }\rightrightarrows v_{0}$.

\item If $v_{\varepsilon }$ is in $L^{2}\left( \Sigma _{\varepsilon }\right)
$ such that
\begin{equation*}
\varepsilon \Vert v_{\varepsilon }\Vert _{L^{2}\left( \Sigma
_{\varepsilon }\right) }^{2}\leq C\text{,}
\end{equation*}%
then one can extract a subsequent still denoted by $\varepsilon $
and a
two-scale limit $v_{0}\in L^{2}((\Omega \times \Sigma )$ such that $%
v_{\varepsilon }\overset{S}{\rightrightarrows }v_{0}$.
\end{enumerate}
\end{proposition}

\subsection{Two-scale limit system}

By virtue of the estimate (\ref{dc}) and the proposition \ref{p1},
there exists $f\in L^{2}\left( \Omega \times Y\right) $ and $g\in
L^{2}\left( \Omega \times \Sigma \right) $ such that, up to a
subsequence, one has
\begin{equation}
\chi \left( \frac{x}{\varepsilon }\right) f_{\varepsilon }\left(
x\right) \rightrightarrows \chi \left( y\right) f\left( x,y\right)
,\ g_{\varepsilon }\left( x\right) \overset{S}{\rightrightarrows
}g\left( x,y\right) . \label{fg}
\end{equation}

Furthermore we have .

\begin{lemma}
\label{l2} \cite{all}, \cite{ngue1}, \cite{adh}. Let
$u_{\varepsilon }$ be the solution of (\ref{wf}). Then there
exists a subsequence still denoted by $\varepsilon $ and two
functions $u\left( x\right) \in H_{0}^{1}\left( \Omega \right) $,
$u_{1}\left( x,y\right) \in L^{2}(\Omega
;H_{\#}^{1}(Y_{s})/\mathbb{R})$ such that
\begin{equation}
\chi \left( \frac{x}{\varepsilon }\right) u_{\varepsilon }\left(
x\right) \rightrightarrows \chi \left( y\right) u\left( x\right)
\text{,} \label{l2_1}
\end{equation}%
\begin{equation}
\chi \left( \frac{x}{\varepsilon }\right) \nabla u_{\varepsilon
}\left( x\right) \rightrightarrows \chi \left( y\right) \left(
\nabla u\left( x\right) +\nabla _{y}u_{1}\left( x,y\right) \right)
\text{.}  \label{l2_2}
\end{equation}%
Moreover we have
\begin{equation}
\lim_{\varepsilon \rightarrow 0}\int_{\Sigma _{\varepsilon
}}\varepsilon u_{\varepsilon }(x)\varphi \left(
x,\frac{x}{\varepsilon }\right) d\sigma _{\varepsilon }\left(
x\right) =\int_{\Omega }\int_{\Sigma }u(x)\varphi (x,y)d\sigma
\left( y\right) dx  \label{l2_3}
\end{equation}%
for every $\varphi \in C(\overline{\Omega };C_{\#}(Y_{s}))$.
\end{lemma}

\begin{lemma}
\label{l3}We have
\begin{equation*}
\lim_{\varepsilon \rightarrow 0}\int_{\Sigma _{\varepsilon
}}u_{\varepsilon }(x)\varphi \left( x\right) \alpha \left(
\frac{x}{\varepsilon }\right) d\sigma _{\varepsilon }\left(
x\right) =\int_{\Omega }\int_{\Sigma }u_{1}(x)\varphi (x)\alpha
\left( y\right) d\sigma \left( y\right) dx
\end{equation*}%
for all $\varphi \left( x\right) \in C\left( \overline{\Omega
}\right) $.
\end{lemma}

\begin{proof}
Define a function $\psi \left( y\right) \in $
$H_{\#}^{1}(Y_{s})/\mathbb{R}$ solution of the problem
\begin{equation}
\left\{
\begin{array}{c}
-\Delta \theta \left( y\right) =0\text{ in }Y_{s}\text{,} \\
\left( \nabla \theta \left( y\right) \right) \cdot \nu \left(
y\right) =\alpha \left( y\right) \text{ on }\Sigma \text{,} \\
y\longmapsto \theta \left( y\right) \text{ Y-periodic.}%
\end{array}%
\right.   \label{pr1}
\end{equation}%
Such a function exists since $\alpha (y)$ satisifies (\ref{as1})
which is the compatibility condition for the solvability of the
problem (\ref{pr1}). Set $\psi \left( y\right) =\nabla \theta $
and consider the function $\psi _{\varepsilon }\left( x\right)
=\psi \left( \frac{x}{\varepsilon }\right) $. Then we have
\begin{eqnarray}
\int_{\Omega _{\varepsilon }}\nabla u_{\varepsilon }(x)\psi
_{\varepsilon }\left( x\right) dx &=&\int_{\Sigma _{\varepsilon
}}u_{\varepsilon }(x)\psi
_{\varepsilon }\left( x\right) \cdot \nu \left( \frac{x}{\varepsilon }%
\right) d\sigma _{\varepsilon }\left( x\right)   \label{le12} \\
&=&\int_{\Sigma _{\varepsilon }}u_{\varepsilon }(x)\alpha \left( \frac{x}{%
\varepsilon }\right) d\sigma _{\varepsilon }\left( x\right) .
\notag
\end{eqnarray}%
Passing to the limit in the left hand side of (\ref{le12}) and
taking into account (\ref{l2_2}) we find
\begin{equation}
\underset{\varepsilon \rightarrow 0}{\lim }\int_{\Sigma
_{\varepsilon }}u_{\varepsilon }(x)\alpha \left(
\frac{x}{\varepsilon }\right) d\sigma _{\varepsilon }\left(
x\right) =\int_{\Omega }\int_{Y}\chi \left( y\right) \left( \nabla
u\left( x\right) +\nabla _{y}u_{1}\left( x,y\right) \right) \psi
\left( y\right) dydx.  \label{le13}
\end{equation}%
Since $u\in H_{0}^{1}\left( \Omega \right) $ we have
\begin{equation*}
\int_{\Omega }\int_{Y}\chi \left( y\right) \nabla u\left( x\right)
\psi \left( y\right) dydx=\left( \int_{\Omega }\nabla u\left(
x\right) dx\right) \int_{Y}\chi \left( y\right) \psi \left(
y\right) dy=0.
\end{equation*}%
Hence the right hand side of (\ref{le13}) becomes
\begin{equation*}
\int_{\Omega }\int_{Y}\chi \left( y\right) \nabla _{y}u_{1}\left(
x,y\right) \psi \left( y\right) dydx.
\end{equation*}%
On the other hand, we have
\begin{eqnarray*}
\int_{\Omega }\int_{Y}\chi \left( y\right) \nabla _{y}u_{1}\left(
x,y\right) \psi \left( y\right) dydx &=&-\int_{\Omega
}\int_{Y_{1}}u_{1}\left( x,y\right) div_{y}\psi \left( y\right)
dydx \\ &&+\int_{\Omega }\int_{\Sigma }u_{1}\left( x,y\right) \psi
\left( y\right) \cdot \nu \left( y\right) d\sigma \left( y\right)
dx \\ &=&\int_{\Omega }\int_{\Sigma }u_{1}\left( x,y\right) \alpha
\left( y\right) d\sigma \left( y\right) dx
\end{eqnarray*}%
which proves the Lemma.
\end{proof}

Now we are able to give the two-scale limit system:

\begin{proposition}
The couple $\left( u,u_{1}\right) \in H_{0}^{1}\left( \Omega
\right) \times L^{2}\left( \Omega ;H_{\#}^{1}\left( Y_{s}\right)
/\mathbb{R}\right) $ is
the solution of the following two-scale homogenized system :%
\begin{gather}
-div_{y}\left( A\left( \nabla u+\nabla _{y}u_{1}\right) \right)
=0\text{\ \ \ in }\Omega \times Y_{s},  \label{tsl1} \\ \left(
A\left( \nabla u+\nabla _{y}u_{1}\right) \cdot \nu \right) +\alpha
u=0\ \ \text{ on }\Omega \times \Sigma \text{,}  \label{tsl2} \\
y\longmapsto u_{1}\text{\ \ \ }Y-\text{periodic,}  \label{tsl3} \\
-div_{x}\left( \int_{Y_{1}}A\left( \nabla u+\nabla
_{y}u_{1}\right) dy\right) +\tilde{\mu}u+\int_{\Sigma }\alpha
u_{1}d\sigma \left( y\right) =F\ \ \ \text{in}\ \Omega \text{,}
\label{tsl4} \\ u=0\text{\ \ \ on }\Gamma   \label{tsl5}
\end{gather}%
where $\tilde{\mu}=\int_{Y_{s}}\mu \left( y\right) dy$ and $%
F(x)=\int_{Y}\chi \left( y\right) f(x,y)dy+\int_{\Sigma }g\left(
x,y\right) d\sigma \left( y\right) .$
\end{proposition}

\begin{proof}
Let $\varphi \left( x\right) \in \mathcal{D}\left( \Omega \right) $ and $%
\varphi _{1}\left( x,y\right) \in \mathcal{D}\left( \Omega
;C_{\#}^{\infty }\left( Y_{s}\right) \right) $. Choosing $v\left(
x\right) =\varphi \left( x\right) +\varepsilon \varphi _{1}\left(
x,\frac{x}{\varepsilon }\right) $
as a test function in problem (\ref{wf}), we have%
\begin{gather}
\int_{\Omega }\nabla u\left( x\right) \chi \left( \frac{x}{\varepsilon }%
\right) ^{t}A(\frac{x}{\varepsilon })\left( \nabla \varphi \left(
x\right) +\varepsilon \nabla _{x}\varphi _{1}\left(
x,\frac{x}{\varepsilon }\right) +\nabla _{y}\varphi _{1}\left(
x,\frac{x}{\varepsilon }\right) \right) dx \notag \\ +\int_{\Omega
_{\varepsilon }}\mu _{\varepsilon }\left( x\right) u\left(
x\right) \left[ \varphi \left( x\right) +\varepsilon \varphi _{1}\left( x,%
\frac{x}{\varepsilon }\right) \right] dx  \notag \\ +\int_{\Sigma
_{\varepsilon }}\alpha _{\varepsilon }\left( x\right) u\left(
x\right) \left[ \varphi \left( x\right) +\varepsilon \varphi _{1}\left( x,%
\frac{x}{\varepsilon }\right) \right] d\sigma _{\varepsilon
}\left( x\right) =  \label{e4} \\
\int_{\Omega }\chi \left( \frac{x}{\varepsilon }\right) f_{\varepsilon }(x)%
\left[ \varphi \left( x\right) +\varepsilon \varphi _{1}\left( x,\frac{x}{%
\varepsilon }\right) \right] dx  \notag \\ \varepsilon
\int_{\Sigma _{\varepsilon }}g_{\varepsilon }\left( x\right)
\left[ \varphi \left( x\right) +\varepsilon \varphi _{1}\left( x,\frac{x}{%
\varepsilon }\right) \right] d\sigma _{\varepsilon }\left(
x\right) .  \notag
\end{gather}%
By virtue of (\ref{l2_2}) the first two terms of the left hand side of (\ref%
{e4}) converges to
\begin{eqnarray}
&&\int_{\Omega }\int_{Y_{s}}[A\left( y\right) \left[ \nabla
u\left( x\right) +\nabla _{y}u_{1}\left( x,y\right) \right] \left[
\nabla \varphi \left( x\right) +\nabla _{y}\varphi _{1}\left(
x,y\right) \right]   \notag \\ &&+\mu \left( y\right) u\left(
x\right) \varphi \left( x\right) ]dydx. \label{r1}
\end{eqnarray}%
Taking into account (\ref{l2_3}), and the lemma \ref{l3}, the
third term of the left hand side of (\ref{e4}) tends to
\begin{equation}
\int_{\Omega }\int_{\Sigma }\alpha \left( y\right) u_{1}\left(
x,y\right) \varphi \left( x\right) d\sigma \left( y\right)
dx+\int_{\Omega }\int_{\Sigma }\alpha \left( y\right) u\left(
x\right) \varphi _{1}\left( x,y\right) d\sigma \left( y\right) dx.
\label{r2}
\end{equation}%
Thanks to (\ref{fg}) the right hand side of (\ref{e4}) converges
to
\begin{equation}
\int_{\Omega }\left[ \int_{Y_{s}}f(x,y)dy+\int_{\Sigma }g\left(
x,y\right) d\sigma \left( y\right) \right] \varphi \left( x\right)
dx=\int_{\Omega }F(x)\varphi \left( x\right) dx.  \label{r3}
\end{equation}

By the density of $\mathcal{D}\left( \Omega \right) $ $\times \in \mathcal{D}%
\left( \Omega ;C_{\#}^{\infty }\left( Y_{s}\right) \right) $ in $%
H_{0}^{1}\left( \Omega \right) \times L^{2}\left( \Omega
;H_{\#}^{1}\left(
Y_{s}\right) /\mathbb{R}\right) $ we get from the limits (\ref{r1})-(\ref{r3}%
) the following two-scale weak formulation system:%
\begin{equation}
\left\{
\begin{array}{c}
\left( u,u_{1}\right) \in H_{0}^{1}\left( \Omega \right) \times
L^{2}\left( \Omega ;H_{\#}^{1}\left( Y_{s}\right)
/\mathbb{R}\right) \text{ is such that } \\ \int_{\Omega
}\int_{Y_{s}}A\left[ \nabla u+\nabla _{y}u_{1}\right] \left[
\nabla v+\nabla _{y}v_{1}\right] dydx+ \\ \tilde{\mu}\int_{\Omega
}uvdx+\int_{\Omega }\int_{\Sigma }\alpha u_{1}vd\sigma \left(
y\right) dx+\int_{\Omega }\int_{\Sigma }\alpha
uv_{1}d\sigma \left( y\right) =\int_{\Omega }Fvdx%
\end{array}%
\right.   \label{e5}
\end{equation}%
for all $\left( v,v_{1}\right) \in H_{0}^{1}\left( \Omega \right)
\times L^{2}\left( \Omega ;H_{\#}^{1}\left( Y_{s}\right)
/\mathbb{R}\right) $.
Integration by parts in (\ref{e5}) \ with respect to $v_{1}$ ($v=0$) gives (%
\ref{tsl1})- (\ref{tsl3}) and with respect to $v$ ($v_{1}=0$) yields (\ref%
{tsl4})-(\ref{tsl5}). The Proposition is now proved.
\end{proof}

Thanks to the linearity of the first equation of (\ref{tsl1}) we
can compute
$u_{1}(x,y)$ in terms of $u\left( x\right) $ as follows:%
\begin{equation}
u_{1}(x,y)=\underset{k=1}{\overset{n}{\sum }}\zeta _{k}\left( y\right) \frac{%
\partial u}{\partial x_{k}}\left( x\right) +\gamma (y)u\left( x\right) +%
\tilde{u}\left( x\right)   \label{rel1}
\end{equation}%
where for each $k$ the function $\zeta _{k}\left( y\right) $
satisfies the
auxiliary problem:%
\begin{eqnarray*}
-div\left( A\left( y\right) \nabla \zeta _{k}\left( y\right)
\right)
&=&div\left( A\left( y\right) e_{k}\right) \text{ in }\Omega \times Y_{s}%
\text{, } \\ A\left( y\right) \nabla \zeta _{k}\left( y\right)
\cdot \nu  &=&-A\left( y\right) e_{k}\cdot \nu \text{ on }\Omega
\times \Sigma \text{,} \\ y &\rightarrow &\zeta _{k}\left(
y\right) \text{ }Y\text{-periodic, }x\in \Omega \text{.}
\end{eqnarray*}%
where $e_{k}=\left( \delta _{ik}\right) _{1\leq i\leq n}$, $\delta
_{ik}$ is the Kr\"{o}necker symbol.

The function $\gamma (y)$ satisfies
\begin{gather*}
-div\left( A\left( y\right) \nabla \gamma \left( y\right) \right)
=0\text{ in }\Omega \times Y_{s}\text{, } \\ A\left( y\right)
\nabla \gamma \left( y\right) \cdot \nu =-\alpha \left( y\right)
\text{ on }\Omega \times \Sigma \text{,} \\ y\rightarrow \gamma
\left( y\right) \text{ }Y\text{-periodic, }x\in \Omega \text{.}
\end{gather*}%
\ \ \newline

Finally, inserting the relation (\ref{rel1}) into the equation
(\ref{tsl4}) yields to the homogenized equation
\begin{equation}
-div\left( A^{\hom }\nabla u(x)\right) +B\cdot \nabla u(x)+\lambda
u\left( x\right) =F\left( x\right)   \label{h1}
\end{equation}%
where $A^{\hom }$ is the matrix with coefficients%
\begin{equation*}
a_{ij}^{\hom }=\underset{k=1}{\overset{n}{\sum
}}\int_{Y_{s}}\left[
a_{ij}\left( y\right) \left( \delta _{kj}+\frac{\partial \zeta _{j}}{%
\partial y_{k}}\left( y\right) \right) \right] dy
\end{equation*}%
$B$ is the vector with components:%
\begin{eqnarray*}
b_{i} &=&\int_{\Sigma }\alpha \left( y\right) \zeta _{i}\left(
y\right)
d\sigma \left( y\right) -\underset{k=1}{\overset{n}{\sum }}%
\int_{Y_{s}}a_{ik}\left( y\right) \frac{\partial \gamma }{\partial y_{k}}%
\left( y\right) dy \\ &=&\int_{\Sigma }\alpha \left( y\right)
\zeta _{i}\left( y\right) d\sigma \left( y\right)
-\int_{Y_{s}}A\left( y\right) e_{i}\nabla \gamma \left( y\right)
dy
\end{eqnarray*}%
$\lambda $ is the real number:%
\begin{eqnarray*}
\lambda  &=&\int_{\Sigma }\alpha \left( y\right) \gamma \left(
y\right) d\sigma \left( y\right) +\tilde{\mu} \\
&=&-\int_{Y_{s}}A\left( y\right) \nabla \gamma \left( y\right)
\nabla \gamma \left( y\right) dy+\tilde{\mu}
\end{eqnarray*}

Thus we have proved the following result

\begin{theorem}
Let $u_{\varepsilon }$ be the solution in $V_{\varepsilon }$ of
the Robin boundary problem (\ref{eq1})-(\ref{eq3}). Then $\chi
_{\varepsilon }\left( x\right) u_{\varepsilon }\left( x\right) $
two-scale converges to $\chi \left( y\right) u\left( x\right) $
where $u\left( x\right) $ is a solution
in $H_{0}^{1}\left( \Omega \right) $ of the homogenized problem:%
\begin{equation}
\left\{
\begin{array}{l}
-div\left( A^{\hom }\nabla u(x)\right) +B\cdot \nabla u(x)+\lambda
u\left( x\right) =F\left( x\right) \text{ in }\Omega \text{,} \\
u=0\text{ on }\Gamma \text{.}%
\end{array}%
\right.   \label{hom}
\end{equation}
\end{theorem}

\begin{remark}
We observe that the limit equation (\ref{hom}) contains an extra
strange
term of order $1$. Namely the convection term $B\cdot \nabla u$. The vector $%
B$ depends closely on the matrix $A$ and the resistivity function
$a$. For
example, if $A$ is symmetric then $B=0$. Indeed%
\begin{equation*}
\int_{\Sigma }\alpha \left( y\right) \zeta _{i}\left( y\right)
d\sigma =-\int_{Y_{s}}A\left( y\right) \nabla \gamma \left(
y\right) \nabla \zeta _{i}\left( y\right) dy=-\int_{Y_{s}}A\left(
y\right) \nabla \zeta _{i}\left( y\right) \nabla \gamma \left(
y\right) dy\text{)}
\end{equation*}%
and since $A$ is symmetric
\begin{equation*}
\int_{\Sigma }\alpha \left( y\right) \zeta _{i}\left( y\right)
d\sigma =\int_{Y_{s}}A\left( y\right) e_{i}\nabla \gamma \left(
y\right) dy.
\end{equation*}

\end{remark}


\begin{thebibliography}{99}
\bibitem{ain} A. Ainouz, Homogenization of a Wentzell type like problems in
elasticity ( In French), Thesis, Algiers, 1997.

\bibitem{ain1} A. Ainouz, Derivation of a double-diffusion model in
poro-elastic media, Proceeding of the First Indo-German Conference
on PDE, Scientific Computing and Optimization in Applications,
September 8-10, 2004, Univ. of Trier, Germany.

\bibitem{all} G. Allaire, Homogenization and two-scale convergence, SIAM J.
Math. Anal., Vol. 23, 6,1482-1518, 1992.

\bibitem{adh} G. Allaire, A. Damlamian, U. Hornung, two-scale convergence on
periodic surfaces and applications, In Proceedings of the
International Conference on Mathematical Modelling of Flow through
Porous Media (May 1995), A. Bourgeat et al. eds., pp.15-25, World
Scientific Pub., Singapore (1996).

\bibitem{bpc} A. G. Belyaev,  A. L. Pyatnitskii and G. A. Chechkin,  Averaging
in a perforated domain with an oscillating third boundary condition.
(Russian) Mat. Sb. 192 (2001), no. 7, 3--20; translation in Sb.
Math. 192 (2001), no. 7-8, 933--949

\bibitem{blp} A. Bensoussan, J. L. Lions and G. Papanicolaou, Asymptotic
analysis for periodic structure, North Holland, Amsterdam, 1978.

\bibitem{ciodon1} D. Cioranescu, P. Donato, Homog\'{e}n\'{e}isation du probl%
\`{e}me de Neumann non homog\`{e}ne dans des ouverts perfor\'{e}s.
Asymptotic Anal. 1 (1988), no. 2, 115--138.

\bibitem{ciodon2} D. Cioranescu, P. Donato, On a Robin problem in
perforated domains. Homogenization and applications to material
sciences (Nice, 1995), 123--135, GAKUTO Internat. Ser. Math. Sci.
Appl., 9, Gakkuto, Tokyo, 1995.

\bibitem{cd} D. Cioranescu, P. Donato, An introduction to homogenization,
Oxford Lectures Series in Mathematics and its Applications 17,
Owford, University Press 1999.

\bibitem{dh} P. Deuflhard, R. Hochmuth, Multiscale analysis of
thermoregulation in the human microvascular system, Math. Meth.
Appli. Sci. 27, pp. 2004 971-989.

\bibitem{lnw} D. Lukkassen, G. Nguetseng and P. Wall, Two-scale Convergence, J.
of Pure and Appl. Math. 2, 1, 35-86, 2002.

\bibitem{nr} M. Neuss-Radu, Some extensions of two-scale convergence. C. R.
Acad. Sci. Paris S\'{e}r. I Math. 322 (1996), no. 9, 899--904.

\bibitem{ngue} G. Nguetseng, A general convergence result for a functional
related to the theory of homogenization, SIAM J. Math. Anal. 20
(1989) 608-623.

\bibitem{ngue1} G. Nguetseng, Asymptotic analysis for a stiff variational
problem arising in mechanics, SIAM J. Math. Anal. Vol. 12 N° 6, pp
1394-1414, 1990.

\bibitem{pas1} S. E. Pastukhova, On the character of the distribution of the
temperature field in a perforated body with a given value on the
outer boundary under heat exchange conditions on the boundary of
the cavities that are in accord with Newton's law. (Russian) Mat.
Sb. 187 (1996), no. 6, 85--96; translation in Sb. Math. 187
(1996), no. 6, 869--880.

\bibitem{san} E. Sanchez-Palencia,  Non-Homogeneous Media and Vibration
Theory, Lecture Notes in Physics 127, 1980.

\bibitem{tar} L. Tartar, Cours Peccot, Coll\`{e}ge de France, 1977.
\end{thebibliography}
\end{document}